\documentclass[graybox,envcountsame,envcountresetchap]{svmult}

\usepackage{graphicx}

\usepackage{amsmath}
\usepackage{amssymb}
\usepackage{color}

\usepackage[latin1]{inputenc}
\usepackage{amsfonts}
\usepackage[english]{babel}
\usepackage[T1]{fontenc}
\usepackage{ae}
\usepackage{url}

\allowdisplaybreaks[4]

\newcommand{\M}[2]{\textnormal{M}\hspace{-0.2em}\left(#1,#2\right)} 

\usepackage{color}
\usepackage{cite}           
\usepackage{mathptmx}       
\usepackage{helvet}         
\usepackage{courier}        
\usepackage{type1cm}        
\usepackage{makeidx}         
\usepackage{graphicx}        
\usepackage{multicol}        
\usepackage[bottom]{footmisc}

\newcommand{\inner}[1]{\langle{#1}\rangle}
\newcommand{\C}{\mathbb{C}}

\newcommand{\Q}{\mathbb{Q}}
\newcommand{\Z}{\mathbb{Z}}
\newcommand{\z}{\zeta}
\newcommand{\al}{\alpha}
\newcommand{\be}{\beta}
\newcommand{\ga}{\gamma}
\renewcommand{\th}{\theta}
\newcommand{\om}{\omega}
\newcommand{\G}{\Gamma}
\newcommand{\new}{\mathrm{new}}
\newcommand{\old}{\mathrm{old}}
\newcommand{\f}{{\mathfrak f}}
\newcommand{\sump}{\sideset{}{'}\sum}
\newcommand{\psmm}[4]{\left(\begin{smallmatrix}{#1}&{#2}\\{#3}&{#4}\end{smallmatrix}\right)}

\DeclareMathOperator{\SL}{SL}
\DeclareMathOperator{\lcm}{lcm}
\newcommand{\T}{{\mathcal T}}
\renewcommand{\E}{{\mathcal E}}
\renewcommand{\H}{\mathfrak H}
\DeclareMathOperator{\Tr}{Tr}

\newcommand{\ov}{\overline}
\newcommand{\leg}[2]{\mbox{$\left(\dfrac{#1}{#2}\right)$}}
\newcommand{\lgs}[2]{\mbox{$\left(\frac{#1}{#2}\right)$}}
\smartqed

\makeindex             

\begin{document}

\title*{Modular Forms\index{Modular Form} in Pari/GP\index{Pari/GP}\\[0.3cm]
{\normalsize\it(Dedicated to Don Zagier for his 65th birthday.)}}
\titlerunning{Modular Forms in Pari/GP}
\author{Karim Belabas and Henri Cohen}
\institute{
Univ. Bordeaux, CNRS, INRIA, IMB, UMR 5251, F-33400 Talence, France.\\
\email{{Henri.Cohen,Karim.Belabas}@math.u-bordeaux.fr}}

\maketitle

\abstract{We give theoretical and practical information on the {\tt Pari/GP}
  modular forms package available since the spring of 2018. Thanks to the use
  of products of two Eisenstein series, this package is the first which can
  compute Fourier expansions at any cusps, evaluate modular forms near the
  real axis, evaluate $L$-functions of non-eigenforms, and compute general
  Petersson scalar products.}

\section{Introduction}

\vspace{1mm}\noindent

Three packages exist which allow computations on classical modular forms:
{\tt Sage}, {\tt Magma}, and {\tt Pari/GP}, the latter being available since
the spring of 2018. The first two packages are based on modular symbols, while
the third is based on trace formulas. This difference is not so important
(although the efficiency of certain computations can vary widely from one
package to another), but at present the {\tt Pari/GP} package is the only one
which is routinely able to perform a number of computations on modular forms
such as expansions at cusps, evaluation near the real axis, evaluation of
$L$-functions of non-eigenforms, computation of general Petersson scalar
products, etc.

The method used for these more advanced commands is based on the one hand on
a theorem of Borisov--Gunnells \cite{Bor-Gun} \cite{Bor-Gun2} stating that
with known exceptions
(which can easily be circumvented) spaces of modular forms are generated by
products of two Eisenstein series, and on the other hand by tedious
computations on the expansions of these Eisenstein series. None of this is
completely original, but it took us several months before obtaining a
satisfactory implementation. In addition, note that we do not need the
Borisov--Gunnells \emph{theorem} (in fact in the beginning we were not even
aware of their work) since we can always check whether products of two
Eisenstein series generate the desired spaces (which they sometimes do not
in weight $2$, but this can be circumvented).

This paper is divided into three parts. In the first part (Sections
\ref{sec:two} to \ref{sec:six}), we describe the theoretical tools used in
the construction of modular form spaces. In the short second part (Section
\ref{sec:seven}), we give some implementation details. In the third somewhat
lengthy part (Sections \ref{sec:eight} and \ref{sec:nine}) we give some sample
commands and results obtained using the package, with emphasis on the advanced
commands not available elsewhere.

\vspace*{4mm}
\noindent
{\bf Acknowledgments.}
We would like to thank B.~Allombert, J.~Bober, A.~Booker, M.~Lee, and
B.~Perrin-Riou for very helpful discussions and help in algorithms and
programming, F.~Brunault and M.~Neururer for Theorem \ref{thmfga}, as well as
K.~Khuri-Makdisi, N.~Mascot, N.~Billerey and E.~Royer.

Last but not least, we thank Don Zagier for his continuous input on this
package and on {\tt Pari/GP} in general. In addition, note that a much smaller
program written 30 years ago by Don, N.~Skoruppa, and the second author
can be considered as an ancestor to the present package.

\section{Construction of Spaces of Integral Weight $k\ge2$}\label{sec:two}

\vspace{1mm}\noindent

\subsection{Introduction}

We decided from the start to restrict to spaces of classical modular forms,
and more precisely to the usual spaces $M_k(\G_0(N),\chi)$ where $\chi$ is
a Dirichlet character modulo $N$ of suitable parity, including $k=1$ and
$k$ half-integral. In addition to this \emph{full} modular form space, we also
want to construct the space of cusp forms $S_k(\G_0(N),\chi)$,
the space of Eisenstein series $\E_k(\G_0(N),\chi)$, 
(so that $M_k(\G_0(N),\chi)=\E_k(\G_0(N),\chi)\oplus S_k(\G_0(N),\chi)$),
the space of newforms $S_k^{\new}(\G_0(N),\chi)$, and the space of oldforms
$S_k^{\old}(\G_0(N),\chi)$ (so that
$S_k(\G_0(N),\chi)=S_k^{\old}(\G_0(N),\chi)\oplus S_k^{\new}(\G_0(N),\chi)$).

Other finite index subgroups of $\SL_2(\Z)$ could be considered, as well as
other subspaces of $M_k(\G_0(N),\chi)$, such as Skoruppa--Zagier's
\emph{certain space} \cite{Sko-Zag}, but to limit the amount of work we
have restricted ourselves to the above. In this section, we assume that
$k \geq 2$ is an integer and defer half-integral weights and weight $1$ to
later sections.

\subsection{Construction of $\E_k(\G_0(N),\chi)$}

The construction of the space of Eisenstein series is easy, and based
on a theorem apparently first published by J.~Weisinger in 1977 \cite{Weis}.
Recall that if $\chi_1$ and $\chi_2$ are two primitive characters modulo
$N_1$ and $N_2$ respectively, we define for $k>2$ the Eisenstein series
$$G_k(\chi_1,\chi_2;\tau)=\sump_{N_1\mid c,\ d}\dfrac{\ov{\chi_1(d)}\chi_2(c/N_1)}{(c\tau+d)^k}\;,$$
and if $k=1$ or $k=2$ one defines $G_k$ by analytic continuation to $s=0$
of the corresponding series where $(c\tau+d)^k$ is multiplied by
$|c\tau+d|^{2s}$ (``Hecke's trick''). Then 
$G_k(\chi_1,\chi_2)$ belongs to $M_k(\G_0(N_1N_2),\chi_1\chi_2)$, except when $k=2$
and $\chi_1$ and $\chi_2$ are trivial characters, in which case there is
a nonanalytic term in $1/\Im(\tau)$.

We introduce the following useful notation: if $\chi$ is any Dirichlet
character, we denote by $\chi_f$ the primitive character equivalent to $\chi$,
and we denote by $1$ the trivial character modulo $1$.
Weisinger's theorem (slightly corrected for $k=2$) is as follows:

\begin{theorem}
  \begin{enumerate}\item[]
  \item For $k\ge3$ or for $k=2$ and $\chi$ a nontrivial character, a
    basis of the space $\E_k(\G_0(N),\chi)$ of Eisenstein series is given by the
    $G_k(\chi_1,\chi_2;m\tau)$, where $(\chi_1,\chi_2)$ ranges over pairs
    of primitive characters as above such that $(\chi_1\chi_2)_f=\chi_f$
    and $N_1N_2\mid N$, and $m$ ranges over all divisors of $N/(N_1N_2)$.
  \item For $k=2$ and $\chi$ a trivial character, a basis of
    $\E_2(\G_0(N))$ is given by the same functions as in (1) except that
    if $(\chi_1,\chi_2)=(1,1)$ we replace $G_k(\chi_1,\chi_2;m\tau)$ by
    $G_2(\chi_1,\chi_2;m\tau)-G_2(\chi_1,\chi_2;\tau)/m$ and exclude $m=1$.
  \item For $k=1$, a basis of $\E_1(\G_0(N),\chi)$ is given by the same
    functions as in (1), except that we restrict to $\chi_1$ being an even
    character.
  \end{enumerate}
\end{theorem}

(Note that since in weight $1$ (and only in weight $1$) the characters
$\chi_1$ and $\chi_2$ play a symmetrical role, we could instead
restrict to $\chi_2$ being an even character.)

Thanks to this theorem it is immediate to construct a basis of
$\E_k(\G_0(N),\chi)$. However, this is not the whole story. Indeed, one can
easily compute the Fourier expansion at infinity of $G_k(\chi_1,\chi_2;\tau)$,
and (after suitable normalization) the coefficients belong to the large
cyclotomic field $\Q(\z_{o_1},\z_{o_2})$,
where $o_i$ denotes the order of the character $\chi_i$ and 
$\z_n$ denotes a primitive $n$th root of unity. It is in fact possible to
obtain a basis whose Fourier coefficients are in the smaller cyclotomic
field $\Q(\z_o)$, where $o$ is the order of $\chi$. For this, we introduce
the following notation:

\begin{definition}\begin{enumerate}
    \item $\Tr_{1,2}$ will denote the \emph{trace map} from
      $\Q(\z_{o_1},\z_{o_2})$ to $\Q(\z_o)$.
    \item We will say that two primitive characters $\chi$ and $\chi'$ modulo
      $N$ are \emph{equivalent} and write $\chi\sim\chi'$ if there exists
      $j$ coprime to the order of $\chi$ such that $\chi'=\chi^j$.
  \end{enumerate}
\end{definition}

\begin{theorem} Let $d_{1,2}=[\Q(\z_{o_1},\z_{o_2}):\Q(\z_o)]$ be the degree
  of the field extension and let $\al_{1,2}$ be such that
  $\Q(\z_{o_1},\z_{o_2})=\Q(\z_o)(\al_{1,2})$. A basis for the space
  $\E_k(\G_0(N),\chi)$ is given by the
  $\Tr_{1,2}(\al_{1,2}^jG_k(\chi_1,\chi_2;m\tau))$
  where $0\le j<d_{1,2}$ and $(\chi_1,\chi_2,m)$ are as in the previous
  theorem (with the suitable modification when $k=2$) except that $\chi_1$
  is only chosen up to equivalence.\end{theorem}

Thus we indeed obtain a basis of the Eisenstein space whose Fourier
expansions have coefficients in the smaller field $\Q(\z_o)$, for instance
in $\Q$ if $\chi$ is trivial or a quadratic character.

\subsection{Construction of $S_k^{\new}(\G_0(N),\chi)$ and of
$S_k(\G_0(N),\chi)$ when $k\geq 2$}

Here, the {\tt Pari/GP} package differs from the others (note that
we do not claim that this is a better choice). First recall the
\emph{Eichler--Selberg trace formula} on $\G_0(N)$. For every $n$ including
those not coprime to
$N$ one defines a Hecke operator $T(n)$ on $M_k(\G_0(N),\chi)$ by the formula
$$T(n)(f)(\tau)=\dfrac{1}{n}\sum_{\substack{ad=n\\\gcd(a,N)=1}}\chi(a)a^k\sum_{b\bmod d}f\left(\dfrac{a\tau+b}{d}\right)\;.$$
Because of the condition $\gcd(a,N)=1$ (which is irrelevant if
$\gcd(n,N)=1$) it is important to note that when $\gcd(n,N)>1$ the
operator $T(n)$ \emph{depends on the level} $N$ of the underlying space,
so should be more properly be denoted $T_N(n)$. Equivalently, if we always
consider $\chi$ as a Dirichlet character modulo $N$, so such that $\chi(a)=0$
when $\gcd(a,N)>1$, we can omit the condition $\gcd(a,N)=1$.

An important formula, due to Selberg and Eichler, gives the \emph{trace}
of $T(n)$ on $S_k(\G_0(N),\chi)$:

\begin{theorem}\label{thmtrace} Let $\chi$ be a Dirichlet character modulo $N$
and let $k\ge2$ be an integer such that $\chi(-1)=(-1)^k$. For all $n\ge1$,
including those not coprime to $N$, we have
$$\Tr_{S_k(\G_0(N),\chi)}(T(n))=A_1-A_2-A_3+A_4\;,$$
where the different contributions $A_i$ are as follows:

$$A_1=n^{k/2-1}\chi(\sqrt{n})\dfrac{k-1}{12}N\prod_{p\mid N}\left(1+\dfrac{1}{p}\right)\;,$$
where it is understood that $\chi(\sqrt{n})=0$ if $n$ is not a square 
(including when $\chi$ is a trivial character).
$$A_2=\kern-5pt\sum_{\substack{t\in\Z\\t^2-4n<0}}\dfrac{\rho^{k-1}-\ov{\rho}^{k-1}}{\rho-\ov{\rho}}\sum_{f^2\mid (t^2-4n)}\dfrac{h((t^2-4n)/f^2)}{w((t^2-4n)/f^2)}\mu(t,\gcd(N,f),n)\;,$$
with
$$\mu(t,g,n)=g\prod_{\substack{p\mid N\\p\nmid N/g}}\left(1+\dfrac{1}{p}\right)\sum_{\substack{x\bmod N\\x^2-tx+n\equiv0\pmod{Ng}}}\chi(x)\;,$$
where $\rho$ and $\ov{\rho}$ are the roots of the polynomial $X^2-tX+n$, in other words, $\rho+\ov{\rho}=t$ and $\rho\ov{\rho}=n$, 
and for $d<0$, $h(d)$ and $w(d)$ are the class number and number of roots of
unity of the quadratic order of discriminant $d$.

$$A_3=\sump_{\substack{d\mid n\\d\le n^{1/2}}}d^{k-1}\kern-10pt\sum_{\substack{c\mid N\\\gcd(c,N/c)\mid\gcd(N/\f(\chi),n/d-d)}}\phi(\gcd(c,N/c))\chi(x_1)\;,$$
where:
\begin{itemize}\item $\sump$ means that the term $d=n^{1/2}$, if present, must
be counted with coefficient $1/2$,
\item $\f(\chi)$ is the conductor of $\chi$,
\item $x_1$ is defined modulo $\lcm(c,N/c)=N/\gcd(c,N/c)$ by the Chinese
remainder congruences $x_1\equiv d\pmod{c}$ and $x_1\equiv n/d\pmod{N/c}$,
\item $\phi$ is Euler's totient function.
\end{itemize}

$A_4=0$ if either $k>2$ or if $k=2$ and $\chi$ is not the trivial character,
and otherwise, if $k=2$ and $\chi$ is trivial then
$$A_4=\sum_{\substack{t\mid n\\\gcd(n/t,N)=1}}t\;.$$
\end{theorem}

We emphasize that in all the above formulas $\chi(x)=0$ if
$\gcd(x,N)>1$, i.e., $\chi$ is always considered as a character modulo $N$.

From this, a nontrivial application of the M\"obius inversion formula
(explained to us by J.~Bober, A.~Booker, and M.~Lee) allows us to compute the
trace on the new space $S_k^{\new}(\G_0(N),\chi)$.
To simplify notation, denote by $\Tr(N,n)$ (resp., $\Tr^{\new}(N,n)$)
the trace of $T_N(n)$ on $S_k(\G_0(N),\chi)$ (resp.,
$S_k^{\new}(\G_0(N),\chi)$). We introduce the following definitions:

\begin{definition}\begin{enumerate}
  \item We define the multiplicative arithmetic function $\be(n)$ on prime
    powers by $\be(p)=-2$, $\be(p^2)=1$, and $\be(p^a)=0$ for $a\ge3$.
  \item For $m\ge1$ we define the multiplicative arithmetic functions
    $\be_m(n)$ on prime powers by $\be_m(p^a)=\be(p^a)$ if $p\nmid m$ and
    $\be_m(p^a)=\mu(p^a)$ if $p\mid m$, where $\mu$ is the usual M{\"o}bius
    function.
  \item An integer $N$ is said to be \emph{squarefull} if for all primes
    $p\mid N$ we have $p^2\mid N$.
\end{enumerate}\end{definition}

\begin{theorem} Let $\chi_N$ be a Dirichlet character modulo $N$ of
conductor $\f\mid N$, and $k\ge 2$ be an integer such that $\chi_N(-1)=(-1)^k$.
Denote as above by $\chi_{\f}$ the primitive character modulo $\f$ equivalent
to $\chi_N$. Finally, write $N=N_1N_2$ with $\gcd(N_1,N_2)=1$, $N_1$
squarefree and $N_2$ squarefull. We have
$$\Tr^{\new}(N,n)=\sum_{\f\mid M\mid N}\kern5pt\sum_{\substack{d\mid\gcd(M/\f,N_1)\\d^2\mid n}}\chi_{\f}(d)d^{k-1}\be_{n/d^2}(N/M)\Tr(M/d,n/d^2).$$
\end{theorem}

The point of this theorem is the following: set
$$\T^{\new}(N)=\sum_{n\ge1}\Tr^{\new}(N,n)q^n\;.$$
Then $\T^{\new}(N)$ is equal to the sum of the normalized eigenforms in
$S_k^{\new}(\G_0(N),\chi)$, hence a simple argument shows that the
$T(n)\T^{\new}(N)$ generate $S_k^{\new}(\G_0(N),\chi)$, so we simply construct
these forms until the dimension of the space they generate is equal to
the dimension of the full new space (equal to $\Tr^{\new}(N,1)$).

\medskip

Once we have obtained a basis for the space $S_k^{\new}(\G_0(N),\chi)$,
it is immediate to obtain a basis of $S_k(\G_0(N),\chi)$ thanks to the
relation
$$S_k(\G_0(N),\chi)=\bigoplus_{\f\mid M\mid N}\bigoplus_{d\mid N/M}B(d)S_k^{\new}(\G_0(M),\chi_{\f})\;,$$
where $B(d)$ is the usual expanding operator $\tau\mapsto d\tau$.

The old space is given by the same formula but restricting to $M<N$:
$$S_k^{\old}(\G_0(N),\chi)=\bigoplus_{\substack{\f\mid M\mid N\\M<N}}\bigoplus_{d\mid N/M}B(d)S_k^{\new}(\G_0(M),\chi_{\f})\;,$$

Note that one could think of using directly the trace formula on the full
cuspidal space (Theorem \ref{thmtrace}), but experiment and complexity
analysis both show that, in addition to being much less canonical, it would
also be less efficient. (Since on the one hand the $A_i$ to be computed are
eventually the same, and on the other hand linear algebra's cost is
superlinear in the dimension, it is more costly to work in a direct sum than in
each subspace independently.)

\section{Construction of Modular Forms of Half-Integral Weight}\label{sec:three}

Recall that modular form spaces of half-integral weight $M_k(\G_0(N),\chi)$
with $k\in1/2+\Z$ are defined only when $4\mid N$ and $\chi$ is an even character.
In weight $1/2$ a beautiful theorem of Serre--Stark asserts that
the space $M_{1/2}(\G_0(N),\chi)$ is spanned by unary theta series, and
the theorem also specifies the cuspidal subspace $S_{1/2}(\G_0(N),\chi)$.
Thus we
consider the construction of the spaces $S_k(\G_0(N),\chi)$ and
$M_k(\G_0(N),\chi)$ when $4\mid N$, $k\ge3/2$ is a half-integer, and $\chi$
is an even character.

Recall the standard theta series
$$\th(\tau)=\sum_{n\in\Z}q^{n^2}=1+2\sum_{n\ge1}q^{n^2}\in M_{1/2}(\G_0(4))\;.$$
Because of the well-known product expansion
$$\th(\tau)=\prod_{n\ge1}(1-q^{2n})(1+q^{2n-1})^2$$
it is clear that $\th$ does not vanish on the
upper half-plane $\H$, and it does not vanish at the cusps $i\infty$ and $0$
of $\G_0(4)$. On the other hand, since the cusp $1/2$ is \emph{irregular},
$\th$ necessarily vanishes at the cusp $1/2$. Applying $\G_0(4)$, we see
that $\th(\tau)=0$ if and only if $\tau$ is a cusp of the form $a/b$
with $\gcd(a,b)=1$ and $b\equiv2\pmod4$.

It follows that $\th(2\tau)=0$ if and only if $\tau$ is a cusp of the form
$a/b$ with $b\equiv4\pmod8$. In particular we see the essential fact that
$\th(\tau)$ and $\th_2(\tau)=\th(2\tau)$ have \emph{no common zeros} in the
completed upper half-plane (we say that they are \emph{coprime forms}).

This allows us to construct the desired modular form spaces of half-integral
weight as follows. Let $k\in\Z+1/2$ and say we want to construct
$M_k(\G_0(N),\chi)$. If $f\in M_k(\G_0(N),\chi)$ then
$f\th\in M_{k+1/2}(\G_0(N),\chi')$, where
$\chi'=\chi$ if $k+1/2\equiv0\pmod2$ otherwise $\chi'=\chi\chi_{-4}$,
where $\chi_{-4}(n)=\leg{-4}{n}$, and similarly
$$f\th_2\in M_{k+1/2}(\G_0(N'),\chi')\supset M_{k+1/2}(\G_0(N),\chi')\;,$$
where $N'=N$ if $8\mid N$, and $N'=2N$ otherwise. By the preceding section we
know how to construct a basis $B$ of $M_{k+1/2}(\G_0(N'),\chi')$.

Now the forms $g_1=f\th$ and $g_2=f\th_2$ which are both in that space
satisfy $g_1\th_2=g_2\th$. This equality can be solved by simple linear
algebra on the basis $B$, and once a basis of $(g_1,g_2)$ is found one
recovers $f$ as the quotient $g_1/\th$ (or $g_2/\th_2$).
Since the level $N'$ is at most twice the initial level $N$, this gives an
efficient method for computing $M_k(\G_0(N),\chi)$. To compute
the cuspidal space $S_k(\G_0(N),\chi)$, simply replace all the $M_{k+1/2}$ by
$S_{k+1/2}$.

Note that we have a check on the correctness of the result by using a
theorem of Oesterl\'e and the second author \cite{Coh-Oes} which gives
the dimensions of $M_k(\G_0(N),\chi)$ and $S_k(\G_0(N),\chi)$ when
$k\in1/2+\Z$.

\medskip

The reader has certainly noticed that we do not speak of the old/new space,
nor of the Eisenstein space. The construction of the old/new space is better
performed in the so-called \emph{Kohnen $+$-space}, and is implemented in the
package, but we do not explain the details here.

It should be possible to construct the Eisenstein space explicitly in
a manner analogous to Weisinger's theorem, but as far as the authors are
aware this has been done only when $N/4$ is squarefree.

\section{Construction of Modular Forms of Weight $1$}\label{sec:four}

Although in principle it is algorithmically just as simple to construct
modular forms of weight $1$ as modular forms of half-integral weight,
it is more difficult to do it efficiently.

A first method which comes to mind is again to use two coprime forms.
In fact, we can again use $\th$ and $\th_2$ as in the previous section,
but to stay in the realm of integral weight forms it is preferable to use
the weight $1$ coprime forms $\th^2$ and $\th_2^2$. This works
exactly in the same way that we used for half-integral weight, but the main
efficiency loss is due to the level: since $\th_2^2\in M_1(\G_0(8),\chi_{-4})$,
the level of $f\th_2^2$ will be $N' = \lcm(8,N)$. In the half-integral case
we always had $4\mid N$, so $N'$ was at most $2N$, but here if
for instance $N$ is odd, $N'=8N$ so we are required to work in
a space of weight $2$ forms and level $8$ times larger, which is
prohibitive since the complexity is at least proportional to the cost of
linear algebra in dimension $N'$.

We can search for other coprime forms. For instance J.~Bober (personal
communication) suggests to use two specific Eisenstein series of weight $1$
and levels $3$ and $4$ respectively. We would then need to work in level
$\lcm(12,N)$, which can be lower than $\lcm(8,N)$ when $3\mid N$ for instance.

\medskip

However, to our knowledge the most efficient method to construct spaces of
modular forms of weight $1$, and the one which is implemented in the
{\tt Pari/GP} package, is the use of Schaeffer's \emph{Hecke stability}
theorem. This theorem essentially states the following: if $V$ is a
\emph{finite-dimensional} vector space of meromorphic modular functions
over $\G_0(N)$ with character $\chi$, and if $V$ is stable by any single
Hecke operator $T(n)$ with $n$ coprime to $N$, then $V$ is in fact a space of
holomorphic modular forms.

Since Eisenstein series of weight $1$ are just as explicit as in higher
weight, it is sufficient to construct the cuspidal space $S_1(\G_0(N),\chi)$,
and to do so we proceed as follows. Let $\E_1(\G_0(N),\ov{\chi})$ be the
space of Eisenstein series of weight $1$ and conjugate character.
Note that if $f\in S_1(\G_0(N),\chi)$ and $E\in\E_1(\G_0(N),\ov{\chi})$
then $fE\in S_2(\G_0(N))$. Consider the (finite dimensional) space
$$W=\bigcap_{E\in B_1(\G_0(N),\ov{\chi})}\dfrac{S_2(\G_0(N))}{E}\;,$$
where $B_1(\G_0(N),\ov{\chi})$ is a basis of $\E_1(\G_0(N),\ov{\chi})$.
It is clear that $S_1(\G_0(N),\chi)\subset W$, and sometimes $W$ is $0$
(for instance if $N\le22$), so we are done. In general this is not the
case, so we apply Schaeffer's theorem. It is easy to show that the
\emph{maximal} stable subspace of $W$ under the action of $T(n)$ (for some
fixed $n$ coprime to $N$) is exactly equal to the desired space
$S_1(\G_0(N),\chi)$.

All the above operations (intersections of spaces and finding maximal stable
subspaces) are elementary linear algebra, but can be extremely expensive
in particular when the values of the character $\chi$ lie in a large
cyclotomic field. Even when $\chi$ has small order, the computations
suffer from intermediate coefficient explosion whereas we expect the
final result to have tiny dimension. We thus use modular algorithms and
perform the computations in various \emph{finite fields} before lifting the
final result.

Note that in the actual implementation we first look for \emph{dihedral
forms}, i.e., forms coming from Hecke Gr\"ossencharacters of quadratic
fields. Once these forms computed, the orders of the possible characters
for the so-called \emph{exotic} forms is much more limited.

\section{Elementary Computations on Modular Forms}\label{sec:elem}\label{sec:five}

Note that since the construction of our modular form spaces ultimately
boils down to the computation of the \emph{trace forms} $\T^{\new}$,
modular forms are always implicitly given by their Fourier expansion at
infinity, which can be unfortunate for some applications.

Nonetheless, a large number of standard operations can be done on modular
forms represented in this way: first, elementary arithmetic operations such as
products, quotients, linear combinations, derivatives, etc., and second,
specifically modular operations such as the action of Hecke operators, of
the expanding operator $B(d)$, Rankin--Cohen brackets, twisting, etc.

Several limitations immediately come to mind:
the action of the Atkin--Lehner operators can be described explicitly only
when the level is squarefree. One can \emph{evaluate} numerically a modular
form by summing its $q$-expansion, but only if $|q|$ is not too close to $1$,
i.e., if $\Im(\tau)$ is not too small (in small levels one can use the action
of $\G_0(N)$ to increase $\Im(\tau)$). One can compute the
Fourier expansion at other cusps than infinity, but only if there exists
an Atkin--Lehner involution sending infinity to that cusp, which will not
always be the case in nonsquarefree level. We will see below how these
limitations can be lifted by the use of the Borisov--Gunnells theorem.

\smallskip

An important operation is \emph{splitting}: once
the new space $S_k^{\new}(\G_0(N),\chi)$ has been constructed, we want to
compute the basis of normalized Hecke eigenforms. This is done by simple
linear algebra after factoring the characteristic polynomials of a sufficient
number of elements of the Hecke algebra. Note that it is in general sufficient
to use the $T(p)$ themselves, but it may happen that one needs more
complicated elements. For instance, to split the space
$S_2^{\new}(\G_0(512))$, no amount of $T(n)$ will be sufficient
(the characteristic polynomials will always have square factors), but
one needs to use in addition for instance the operator $T(3)+T(5)$.

\smallskip

Among the other elementary computations, note that modular forms can arise
naturally from several different sources: forms associated to elliptic curves
defined over $\Q$, forms coming from theta functions of lattices (possibly
with a spherical polynomial), eta quotients, or forms coming
from natural $L$-functions whose gamma factor is
$\G_{\C}(s)=2(2\pi)^{-s}\G(s)$.

\section{Advanced Computations on Modular Forms}\label{sec:six}

We now come to the more advanced functions of the package, which are up to now
not available elsewhere. In view of the limitations above, the basic stumbling
block is the computation of the Fourier expansions of a modular forms
at cusps other than $i\infty$.

For $f\in M_k(\G_0(N),\chi)$, we want more generally to be able to compute
the Fourier expansion of $f|_k\ga$ for any $\ga\in\G$ (in fact it is trivial
to generalize the construction to any $\ga\in M_2^+(\Q)$, and this is done
in the package but will not be explained here). We can assume that $k$ is an
integer: indeed the expansion of $\theta|_k\ga$ is known and we can compute
the expansion of $f\theta$ when $f$ has half-integral weight.
We recall that if
$\ga=\psmm{a}{b}{c}{d}\in\G$ then
$$f|_k\ga(\tau)=(c\tau+d)^{-k}f\left(\dfrac{a\tau+b}{c\tau+d}\right)\;.$$
First, what precisely do we mean by ``Fourier expansion'' ? It is easy to
show that there exists an integer $w\le N$ and a rational number
$\al\in[0,1[\cap\Q$ such that\footnote{$[0,1[$ is a much more sensible notation than $[0,1)$, and $]0,1[$ than $(0,1)$ which can mean so many things.}
$$f|_k\ga(\tau)=q^{\al}\sum_{n\ge0}a_{\ga}(n)q^{n/w}\;,$$
where we always use the convention that $q^x=e^{2\pi ix\tau}$ for $x\in\Q$.
More precisely one can always choose $w=N/\gcd(N,c^2)$, the \emph{width} of
the cusp $a/c=\ga(i\infty)$, and $\al$ is the unique number in $[0,1[$
such that $\chi(1+acw)=e^{2\pi i\al}$. Note in passing that the denominator
of $\al$ divides $\gcd(N,c^2)/\gcd(N,c)$, and by definition that $\al$ is
nonzero if and only if the cusp $a/c$ is so-called \emph{irregular} for the
space $M_k(\G_0(N),\chi)$.

\medskip

The basic idea is as follows. Consider the space generated by products of
two Eisenstein series $E_1$ and $E_2$, chosen of course such that
$E_1E_2\in M_k(\G_0(N),\chi)$, including the trivial Eisenstein series $1$
of weight $0$. Experiments show that usually this space is the whole of
$M_k(\G_0(N),\chi)$, and in fact the only exceptions seem to be in weight
$k=2$. In fact it is a theorem of Borisov--Gunnells \cite{Bor-Gun}
\cite{Bor-Gun2} that
this observation is indeed true, and that the exceptions occur only in
weight $2$ when there exists an eigenform $f$ such that $L(f,1)=0$,
for instance modular forms attached to elliptic curves of positive rank
(so $N=37$ is the smallest level for which there is an exception).

Assume for the moment that we are in a space generated by products of two
Eisenstein series. Since these are so completely explicit, it is possible by
a tedious computation to obtain the expansions of $E|\ga$, hence
of all our forms. In the special case (occurring only in weight $2$) where
the space is not generated by products of two Eisenstein series, we simply
multiply by some known Eisenstein series $E$ so as to be in larger weight, do
the computation there, and finally divide by the expansion of $E|\ga$.
This is what we do in any case in weight $1$ and in half-integral weight.

This may sound straightforward, but as mentioned at the beginning, it
required several months of work first to obtain the correct formulas, and
second to write a reasonably efficient implementation. In fact, we had to
make a choice. With our current choice of Eisenstein series (which may change
in the future), the coefficients of the Eisenstein series $E|\ga$ lie in
a very large cyclotomic field, at worst $\Q(\z_{\lcm(N,\phi(N))})$. It is
possible that we can reduce this considerably, but for now we do not know how
to do this at least in a systematic way. Handling such large algebraic
objects is extremely costly, so we chose instead to work with approximate
complex numbers with sufficient accuracy, and if desired to recognize the
algebraic coefficients at the end using the LLL algorithm. This is reflected
in the command that we will explain below to compute these expansions.

\smallskip

Note however the following theorem, communicated to us by
F.~Brunault and M.~Neururer whom we heartily thank:

\begin{theorem}\label{thmfga} Let $f\in M_k(\G_0(N),\chi)$, denote by
  $M\mid N$ the
  conductor of $\chi$, and assume that the coefficients of the Fourier
  expansion of $f$ at infinity all lie in a number field $K$.
  Then if $\ga=\psmm{A}{B}{C}{D}\in\G$ the Fourier coefficients $a_{\ga}(n)$
  of $f|_k\ga$ lie in the field $K(\z_u)$, where
  $u=\lcm(N/\gcd(N,CD),M/\gcd(M,BC))$.
\end{theorem}
Note that the Fourier coefficients of $f|_k\ga$ can live in a smaller number
field, or even in a smaller cyclotomic extension than that predicted by the
theorem.

\medskip

Once solved the problem of computing expansions of $f|_k\ga$, essentially all
of the limitations mentioned in Section \ref{sec:elem} disappear: we can
evaluate a modular form even very near the real axis (or at cusps), we can
compute the action of the Atkin--Lehner operators in nonsquarefree level, we
can compute general period integrals involving modular forms and in particular
\emph{modular symbols} such as
$$\int_a^b(X-\tau)^{k-2}f(\tau)\,d\tau\;$$
when $k \ge 2$ is integral,
we can numerically evaluate $L$-functions of modular forms which are not
necessarily eigenforms at an arbitrary $s\in\C$, and we can
compute general \emph{Petersson scalar products} thanks to the following
theorem similar to Haberland's:

\begin{theorem}\label{thmpet} Let $k\ge2$ be an integer, let $f$ and $g$
  in $S_k(\G_0(N),\chi)$ be two cusp forms, let
  $$\G=\bigsqcup_{j=1}^r\G_0(N)\ga_j$$ be a right
  coset decomposition, and define $f_j=f|_k\ga_j$ and $g_j=g|_k\ga_j$.
  Finally, for any function $h$ and $a$ and $b$ in the completed upper
  half-plane set
  $$I_n(a,b,h)=\int_a^b \tau^nh(\tau)\,d\tau\;,$$
  the integral being taken along a geodesic arc from $a$ to $b$. Then
  we have
  $$6r(2i)^{k-1}\inner{f,g}_{\G_0(N)}=\sum_{j=1}^r\sum_{n=0}^{k-2}
    (-1)^n\binom{k-2}{n}I_{k-2-n}(0,i\infty,f_j)\ov{I_n(-1,1,g_j)}\;.$$
\end{theorem}

More generally, when at each cusp at least one of the two forms vanishes,
there exists a similar formula which we do not give here.

On the other hand this theorem cannot be applied in weight
$k=1$ or when $k\in1/2+\Z$. However, recent work by D.~Collins \cite{Col}
using a formula of P.~Nelson \cite{Nel} allows to compute Petersson products
in these cases, although less efficiently in general.

\section{Implementation Issues}\label{sec:seven}

Modular form spaces can be represented by a basis, in echelon form or not,
together with suitable linear algebra precomputations allowing fast
recognition of elements of the space (recall the notion of \emph{Sturm bound},
which tells us that if the Fourier coefficients at infinity of two modular
forms belonging to the same space are equal up to some effective bound,
the forms are identical).

The main problem is the representation of modular forms themselves. Because
of our choice of using trace formulas, we must in some way represent the
forms by their Fourier expansion at infinity. Since we do not want to specify
in advance the number of coefficients that we want, a first approach would be
to say that a modular form is a program which, given some $L$, outputs the
Fourier coefficients up to $q^L$. This would be quite
inefficient because of the action of Hecke operators: for instance,
let $p$ be a prime not dividing the level $N$, and assume that
$f(\tau)=\sum_{n\ge0}a(n)q^n$. Then $(T(p)f)(\tau)=\sum_{m\ge0}b(m)q^m$
with $b(m)=a(pn)+p^{k-1}\chi(p)a(n/p)$, where the last term occurs only if
$p\mid n$. Thus if we want $L$ Fourier coefficients of $T(p)f$ we need on
the one hand the coefficients $a(m)$ for $m\le L/p$, but also the coefficients
$a(m)$ for $p\mid m$ and $m\le pL$. All the other Fourier coefficients
$a(m)$ for $m\le pL$ with $p\nmid m$ are not needed, so it would be a waste
to compute them if it can be avoided.

Thus we modify our first approach, and we say that a modular form is a
program which, given some $L$ and step size $d$, outputs the Fourier
coefficients of $q^m$ for $m\le dL$ and $d\mid m$ (hence all the coefficients
if $d=1$). Such a program will define a modular form in our package.

Note that such representation may look like magic: for instance in a
{\tt GP} session, type {\tt D=mfDelta()} which creates the Ramanujan
$\Delta$ function. The output is only one and a half line long and contains
mostly trivial information. Nonetheless, this information is sufficient
to compute the Fourier expansion to millions of terms if desired, using
the fundamental function {\tt mfcoefs(D,n)} since internally the small
information calls a much more sophisticated program which computes the
expansion (note that in this specific case it is faster to compute
the expansion directly using the product formula for the delta function
than to use the modular form package).

\medskip

Additional implementation comments: since the trace formula involves
computing the class numbers $h(D)$ for $D<0$, we use a \emph{cache} method
for those: we first precompute a reasonable number of such class numbers,
then if it is not sufficient, we precompute again a larger number,
and so on. We do similar caching for the factorizations and divisors
of integers.

We also need to represent Dirichlet characters $\chi$. Since the most
frequent are (trivial or) quadratic characters modulo $D$, such a character
will simply be represented by the number $D$ (so $D=1$ or omitted completely
when $\chi$ is the trivial character). More general characters can be
represented in several {\tt Pari/GP} compatible ways, but the preferred
way is to use the \emph{Conrey numbering} which we do not explain here,
so that a general Dirichlet character modulo $N$, primitive or not, is
represented by {\tt Mod(a,N)}, where $\gcd(a,N)=1$.

\section{Pari/GP Commands}\label{sec:eight}

\subsection{Commands involving only Modular Forms}

As already mentioned above, the first basic command is {\tt mfcoefs(f,n)}
which gives the vector of Fourier coefficients $[a(0),a(1),\dotsc,a(n)]$.
We have chosen to give a vector and not a series first because it is more
compact, and second because the series variable (in principle $q$) could
conflict with other user variables. Of course nothing prevents the user
from defining his own function
\begin{verbatim}
? mfser(f,n) = Ser(mfcoefs(f,n),'q);
\end{verbatim}
if the variable $q$ is ok\footnote{A technicality which explains why
representing forms as series with this additional variable is awkward: the
variable $q$ must have higher priority than $t$ otherwise some of the
examples below will fail. A definition which would work in all cases is
{\tt  mfser(f,n) = Ser(mfcoefs(f,n), varhigher("q", 't));}
}. For simplicity we will use this user-defined function in the examples.

\begin{verbatim}
? mfser(mfDelta(), 8)
% = q - 24*q^2 + 252*q^3 - 1472*q^4 + 4830*q^5 
    - 6048*q^6 - 16744*q^7 + 84480*q^8 + O(q^9)
\end{verbatim}

  There are a number of predefined modular forms: in addition to Ramanujan
  $\Delta$ function,
  we have {\tt mfEk(k)}, the normalized Eisenstein series for the full
  modular group $E_k$, more generally {\tt mfeisenstein(k,chi1,chi2)}
  for general Eisenstein series, {\tt mfEH(k)} for Eisenstein series over
  $\G_0(4)$ in half-integral weight, {\tt mfTheta(chi)}, the unary theta
  function associated to the Dirichlet character {\tt chi} (if {\tt chi} is
  omitted, the standard theta series), as well as modular forms coming
  from preexisting mathematical objects such as {\tt mffrometaquo}, eta
  quotients, {\tt mffromell}, modular cusp form of weight $2$ associated to
  an elliptic curve over $\Q$, {\tt mffromqf}, modular form associated to a quadratic
  form, with an optional spherical polynomial, and {\tt mffromlfun}, modular
  form associated to an $L$-function having factor at infinity equal to
  $\G_{\C}(s)=2\cdot(2\pi)^{-s}\G(s)$.

  All the standard arithmetic operations are implemented. For instance,
  the modular form $F=E_4(\Delta^2+E_{24})$ is obtained by the following
  commands:

\begin{verbatim}
? E4 = mfEk(4); E24 = mfEk(24); D = mfDelta();
? F = mfmul(E4, mflinear([mfpow(D,2), E24],[1,1]));
\end{verbatim}

Note that there is no {\tt mfadd}, {\tt mfsub}, or {\tt mfscalmul} functions
since they can be emulated by {\tt mflinear}, which creates arbitrary
linear combinations of forms:

\begin{verbatim}
? mfadd(F,G) = mflinear([F,G],[1,1]);
? mfsub(F,G) = mflinear([F,G],[1,-1]);
? mfscalmul(F,z) = mflinear([F],[z]);
\end{verbatim}

Note also that the internal representation is an expression tree in
\emph{direct Polish} notation.
In fact, since after a number of operations you may have forgotten what
your modular form is, there is a function {\tt mfdescribe} which essentially
outputs this representation. For instance, applying to our above example:

\begin{verbatim}
? mfdescribe(F)
% = "MUL(E_4, LIN([POW(DELTA, 2), E_24], [1, 1]))"
? mfparams(F)
% = [1, 28, 1, y]
\end{verbatim}

  The {\tt mfdescribe} command is not to be confused with the {\tt mfparams}
  command which gives a short description of modularity and arithmetic
  properties of the form. The above command shows that
  $F\in M_{28}(\G_0(1),1)=M_{28}(\G)$, and that the number field
  generated by the Fourier coefficients of $F$ is $\Q[y]/(y) = \Q$.

  Other modular operations are available which work directly on
  forms, such as {\tt mfbd} (expansion $\tau\mapsto d\tau$), {\tt mfderivE2}
  (Serre derivative), {\tt mfbracket} (Rankin--Cohen bracket), etc.
  But most operations need an underlying modular form \emph{space},
  even the Hecke operators, as we have seen above.

\subsection{Commands on Modular Form Spaces}

The basic command which creates a basis of a modular form space is
{\tt mfinit}, analogous to the other {\tt Pari/GP} commands such as
{\tt nfinit}, {\tt bnfinit}, {\tt ellinit}, etc. This command takes two
parameters. The first is a vector {\tt [N,k,CHI]} (or simply {\tt [N,k]}
if {\tt CHI} is trivial), where $N$ is the level, $k$ the weight, and
{\tt CHI} the character in the format briefly described above.
The second parameter is a code specifying which space we want
($0$ for the new space, $1$ for the cuspidal space, or omitted for the
full space for instance). A simpler command with exactly the
same parameters is {\tt mfdim} which gives the dimension.

\begin{verbatim}
? mf = mfinit([26,2],0);
? L = mfbasis(mf); vector(#L,i,mfser(L[i],10))
\end{verbatim}

The first command creates the space $S_2^{\new}(\G_0(26))$ (no character),
and the second commands gives the $q$-expansions of the basis elements:

\begin{verbatim}
% = [2*q - 2*q^3 + 2*q^4 - 4*q^5 -...,\
    -2*q - 4*q^2 + 10*q^3 - 2*q^4 +...]
\end{verbatim}
These are of course not the eigenforms. To obtain the latter:
\begin{verbatim}
? LE = mfeigenbasis(mf); vector(#LE,i,mfser(LE[i],8))
% = [q - q^2 + q^3 + q^4 - 3*q^5 -...,\
     q + q^2 - 3*q^3 + q^4 - q^5 - ...]
\end{verbatim}
Note that eigenforms can be defined over a relative extension of $\Q(\chi)$:

\begin{verbatim}
? mf = mfinit([23,2],0); LE = mfeigenbasis(mf);
? vector(#LE,i,mfser(LE[i],10))
% = [Mod(1, y^2+y-1)*q + Mod(y, y^2+y-1)*q^2 + ...]
\end{verbatim}
There are two ways to better see them:
  
\begin{verbatim}
? mffields(mf)
% = [y^2 + y - 1]
? vector(#LE,i,lift(mfser(LE[i],10)))
% = [q + y*q^2 + (-2*y-1)*q^3 + (-y-1)*q^4 + ...]
? f = LE[1]; mfembed(f, mfcoefs(f,10))
% = [[0, 1,  0.618033988..., -2.236067977..., ...],\
     [0, 1, -1.618033988...,  2.236067977..., ...]]
\end{verbatim}

The first command gives the number fields over which the eigenforms are
defined. Here there is only one eigenform and only one field $\Q[y]/(y^2+y-1)$.
The second command ``lifts'' the coefficients to $\Q[y]$,
so the result is much more legible. The third command ``embeds'' the
eigenform in all possible ways in $\C$: indeed, even though (in the present
example) there is only one \emph{formal} eigenform, the space is of dimension
two so there are \emph{two} eigenforms, given numerically as the last result.

Unavoidably, when nonquadratic characters occur, we can obtain even more
complicated output. We have chosen to represent formal values of a character
(which are in a cyclotomic field) with the variable letter ``{\tt t}'',
but it must be understood that contrary to eigenforms this corresponds to
a single canonical embedding: for instance
\begin{verbatim} Mod(t, t^4 + t^3 + t^2 + t + 1)\end{verbatim}
means in fact $e^{2\pi i/5}$, not some other fifth root of unity.

\begin{verbatim}
? mf = mfinit([15,3,Mod(2,5)], 0); mffields(mf)
% = [y^2 + Mod(-3*t, t^2 + 1)]
? F = mfser(mfeigenbasis(mf)[1], 10); f = liftall(F)
% = q + (y-t-1)*q^2 + t*y*q^3 + ((-2*t-2)*y+t)*q^4...
\end{verbatim}

The first command shows that the eigenforms will have coefficients in
a quadratic extension of a quadratic extension, hence in the quartic field
$\Q[y,t]/(t^2+1, y^2-3t)$.
The second command lifts the $q$-expansion to $\Q[y,t]$.
To obtain the expansion over the quartic field, which is isomorphic to
$\Q[y]/(y^4+9)$, we write
\begin{verbatim}
? [T,a] = rnfequation(t^2+1, y^2-3*t,1);
? T
% = y^4 + 9
? lift(subst(f,t,a))
% = q + (-1/3*y^2+y-1)*q^2 + 1/3*y^3*q^3 + ...
\end{verbatim}
Note that the variable $y$ actually stands for the same algebraic number in
this block and the previous one but {\tt rnfequation} does not guarantee this
in general.
  
\subsection{Miscellaneous Commands}

\begin{verbatim}
? mf = mfinit([96,4],0); M = mfheckemat(mf,5)
% = 
[0    0    64   0   0   -84]
[0    4     0  36   0     0]
[1    0 -24/5   0   0 294/5]
[0    2     0 -12 -20     0]
[0 -1/2     0   1   6     0]
[0    0   6/5   0   0  14/5]

? factor(charpoly(M))
% = [x - 10, 2; x - 2, 2; x + 14, 2]
? M = mfatkininit(mf,3)[2] \\ Atkin-Lehner W_3
% = 
[   0 -3    0    0  -24   0]
[-1/3  0 -4/3    0    0 -12]
[   0  0    0 -9/5 -6/5   0]
[   0  0 -2/3    0    0  -1]
[   0  0  1/6    0    0 3/2]
[   0  0    0  1/5  4/5   0]

? factor(charpoly(M))
% = [x - 1, 3; x + 1, 3]
\end{verbatim}
(all outputs edited for clarity). Self-explanatory. Note that the basis
we compute for our modular form spaces is essentially random and does 
not guarantee that the matrices attached to Hecke or Atkin-Lehner operators
have integral coefficients. They will in general have coefficients in
$\Q(\chi)$ (up to normalizing Gauss sums in the case of $W_Q$); of course,
their characteristic polynomials have coefficients in $\Z[\chi]$.

\begin{verbatim}
? T = mfTheta();
? mf = mfinit([4,2]); mftobasis(mf, mfpow(T, 4))
% = [0, 8]~
? mf = mfinit([4,5,-4]); mftobasis(mf, mfpow(T, 10))
% = [64/5, 4/5, 32/5]~
\end{verbatim}

  Since in both cases the basis of {\tt mf} can be given explicitly, this
  gives explicit formulas for the number of representations of an integer
  as a sum of $r$ squares for $r=4$ and $r=10$ respectively (this can be
  done for $1\le r\le 8$ and $r=10$).

\begin{verbatim}
? B = mfbasis([4,3,-4],3); \\ 3: Eisenstein space
? [mfser(E,10) | E <- B]
% = [q + 4*q^2 + 8*q^3 + 16*q^4 + 26*q^5 + ...\
     -1/4 + q + q^2 - 8*q^3 + q^4 + 26*q^5 + ...]
\end{verbatim}

Modular forms can be evaluated numerically (even when the imaginary part is
very small, see below), as well as their $L$-functions:

\begin{verbatim}
? E4 = mfEk(4); mf = mfinit(E4); mfeval(mf, E4, I)
% = 1.455762892268709322462422003598869...
? 3*gamma(1/4)^8/(2*Pi)^6
% = 1.455762892268709322462422003598869...
\end{verbatim}
  This equality is a consequence of the theory of
  \emph{complex multiplication}, and in particular of the
  \emph{Lerch, Chowla--Selberg} formula.

\begin{verbatim}
? D = mfDelta(); mf = mfinit(D,1); L = lfunmf(mf,D);
? lfunmfspec(L)
% = [[1, 25/48, 5/12, 25/48, 1],\
     [1620/691, 1, 9/14, 9/14, 1, 1620/691],\
      0.0074154209298961305890064277459002287248,\
      0.0050835121083932868604942901374387473226]
\end{verbatim}

The command {\tt lfunmf} creates the $L$-function attached to $\Delta$,
and {\tt lfunmfspec} gives the corresponding \emph{special values}
in the interval $[1,11]$, which are rational numbers times two periods
$\om^+$ (for the odd integers) and $\om^-$ (for the even integers).

{\tt L} is now an $L$-function in the sense of the $L$-function package of
{\tt Pari/GP}, and can be handled as such. For instance:

\begin{verbatim}
? LF=lfuninit(L,[50]); ploth(t=0,50,lfunhardy(LF,t));
\end{verbatim}
outputs in 50 ms the plot of the Hardy function associated to $\Delta$ on
the critical line $\Re(s)=6$ from height $0$ to $50$:

\includegraphics[width=0.85\textwidth]{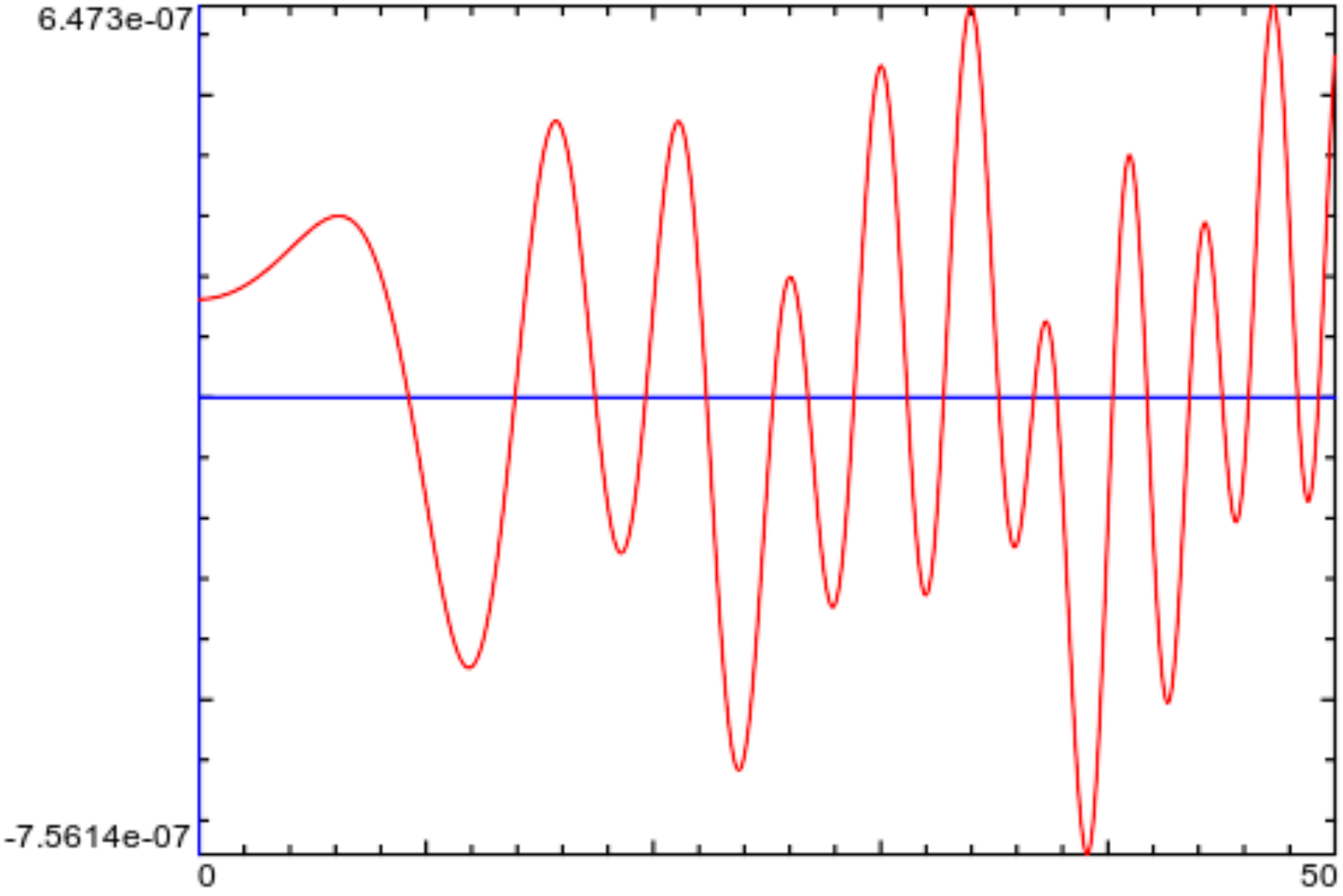}

Similarly, we can compute zeros:

\begin{verbatim}
? lfunzeros(LF,20)
% = [ 9.2223793999211025222437671927434781355,\
     13.907549861392134406446681328770219492,\
     17.442776978234473313551525137127262719,\
     19.656513141954961000127281756321302802]
\end{verbatim}

Note the nontrivial fact that {\tt lfunmf} is completely general and computes
the $L$-function attached to \emph{any} modular form, eigenform or not
(although the computation is indeed more efficient when the function
detects an eigenform): this makes use of the ``advanced'' features that we
will explain below.

\smallskip

A very useful command is {\tt mfeigensearch}, which searches
for \emph{rational} eigenforms (hence with trivial or quadratic character)
in a given range.

\begin{verbatim}
? B = mfeigensearch([[1..60],2],[[2,-1],[3,-3]]);
? apply(mfparams,B)
% [[53, 2, 1, y], [58, 2, 1, y]]
? apply(x->mfser(x,10),B)
% = [q - q^2 - 3*q^3 - q^4 + 3*q^6 +...\
     q - q^2 - 3*q^3 + q^4 - 3*q^5 + 3*q^6 +...]
\end{verbatim}

The first command asks for all rational eigenforms of levels between $1$ and
$60$, of weight $2$ such that $a(2)=-1$ and $a(3)=-3$. The {\tt mfparams}
command shows that there are two such forms, one in level $53$ the other in
level $58$. The last command gives the beginning of their Fourier expansions,
which of course agree up to the coefficient of $q^3$.

Note that the functions are returned as a black box allowing to compute
an arbitrary number of Fourier coefficients, which need not
be specified in advance. We could just as well have written
{\tt mfser(x,1000)} if we wanted $1000$ coefficients.

There exists also the more straightforward {\tt mfsearch} command which
simply searches for a rational modular form with given initial coefficients:

\begin{verbatim}
? B=mfsearch([[1..30],3],[0,1,2,3,4,5,6,7,8],1);
? apply(mfparams,B)
% = [[30, 3, -3, y], [30, 3, -15, y]]
? apply(x->mfser(x,10),B)
% = [q + 2*q^2 + 3*q^3 + 4*q^4 + 5*q^5 + 6*q^6\
       + 7*q^7 + 8*q^8 - 14*q^9 - 30*q^10 + O(q^11),\
     q + 2*q^2 + 3*q^3 + 4*q^4 + 5*q^5 + 6*q^6\
       + 7*q^7 + 8*q^8 - 21*q^9 - 50*q^10 + O(q^11)]
\end{verbatim}

This tells us that there exist exactly two forms of weight $3$ and level
$N\le30$ in the cuspidal space (code $1$) whose Fourier expansion begins by
$q+2q^2+\cdots+8q^8+O(q^9)$; the last command shows their Fourier expansion
up to $q^{10}$.

\subsection{Weight $1$ Examples}

Almost all of the commands given up to now (with the exception of {\tt lfunmf}
for non eigenforms and {\tt mfatkininit} for non squarefree levels) are
direct (although sometimes complicated) applications of the
trace formula and linear algebra over cyclotomic fields. We now come to more
advanced aspects of the package.

As already mentioned, constructing modular forms of weight $1$ is more
difficult than in higher weight, but they are fully implemented in the package.

\begin{verbatim}
? mfdim([148,1,0], 1)
% = [[4, Mod(105, 148), 1, 0],
     [6, Mod(63, 148), 1, 1],
     [18, Mod(127, 148), 1, 1]]
\end{verbatim}
This command uses the \emph{joker} character $0$ (which is available for all
weights but especially useful in weight $1$): it asks to output information
about $S_1(\G_0(148,\chi))$ for all Galois equivalence classes of characters,
but only for nonzero spaces. Here it gives us the Conrey labels of three
characters modulo $148$, such as {\tt Mod(105, 148)}, of respective orders
$4$, $6$ and $18$. The other two integers are of course also important and
give the dimension of the space and the dimension of the subspace generated
by the dihedral forms. Let us look at the first space, which contains an
exotic (non-dihedral) form:
\begin{verbatim}
? mf = mfinit([148,1,Mod(105,148)], 0);
? f = mfeigenbasis(mf)[1];
? mfser(f,12)
% = Mod(1,t^2+1)*q + Mod(-t,t^2+1)*q^3
    + Mod(-1,t^2+1)*q^7 + Mod(t,t^2+1)*q^11 + O(q^13)
? mfgaloistype(mf)
% = [-24]
\end{verbatim}

This tells that the projective image of the Galois representation associated
to the (unique) eigenform in {\tt mf} is isomorphic to $S_4$, so is ``exotic''.
This is the lowest possible level for which it occurs (the smallest
exotic $A_4$ is in level $124$, already found long ago by J.~Tate,
and the smallest exotic $A_5$ is in level $633$ found only a few years ago
by K.~Buzzard and A.~Lauder:

\begin{verbatim}
? mfgaloistype([633, 1, Mod(107,633)])
% = [10, -60]
\end{verbatim}
(The first eigenform is of dihedral type $D_5$, the second has
exotic type $A_5$.) Note that this computation only requires 5 seconds.

\begin{verbatim}
? mfgaloistype([2083,1,-2083])
% = [14, -60]
\end{verbatim}
  This answers an old question of Serre who conjectured the existence of
  exotic $A_5$ forms of prime level $p\equiv3\pmod4$ with quadratic
  character $\lgs{-p}{n}$: $p=2083$ is the smallest such prime.

Typical nonexotic examples:

\begin{verbatim}
? mfgaloistype([239,1,-239])
% = [6, 10, 30]
\end{verbatim}
Three eigenforms with projective image isomorphic to $D_3$, $D_5$, and
$D_{15}$.

\subsection{Half-Integral Weight Examples}

\begin{verbatim}
? mf = mfinit([12,5/2]); B = mfbasis(mf);
? for(j=1,#B,print(mfser(B[j],8)))
1 + 12*q^5 + 30*q^8 + O(q^9)
q - 8*q^5 + 14*q^6 + 28*q^7 - 20*q^8 + O(q^9)
q^2 + 8*q^5 - q^6 - 10*q^7 + 18*q^8 + O(q^9)
q^3 - 4*q^5 + 4*q^6 + 10*q^7 - 10*q^8 + O(q^9)
q^4 + 2*q^5 - 2*q^6 - 4*q^7 + 5*q^8 + O(q^9)
\end{verbatim}

\begin{verbatim}
? f = B[1]; [mf2,F] = mfshimura(mf,f,5); mfser(F,8)
% = -3/5 + 12*q + 108*q^2 + 132*q^3 + 876*q^4 ...
? mfparams(mf2)
% = [6, 4, 1, 4]
\end{verbatim}
This returns the Shimura lift of $f$ of weight $2k-1 = 4$ in $M_4(\G_0(6),
\lgs{4}{\cdot})$. The Kohnen $+$-space as well as the known bijections
between it and spaces of integral weight are implemented, as well as the new
space and the eigenforms in the Kohnen space. We refer the reader to the
manual for details.

\section{Advanced Examples}\label{sec:nine}

\subsection{Expansions of $F|_k\ga$ and Applications}

As mentioned at the beginning, we now give a number of examples which we
believe are not possible (at least in general) with other packages.
The basic function which allows all the remaining advanced examples to
work is the computation of the Fourier expansion of $F|_k\ga$ for any
$\ga\in\G$. We begin with a simple example:

\begin{verbatim}
? mf = mfinit([32,4],0); F = mfbasis(mf)[1];
? mfser(F,10)
% = 3*q + 2*q^5 + 47*q^9 + O(q^11)
? g = [1,0;2,1];
? Ser(mfslashexpansion(mf, F, g, 6, 1, &params), q)
% = Mod(-1/64*t, t^8 + 1)*q
  + Mod(-1/4*t^3, t^8 + 1)*q^3
  + Mod(-11/32*t^5, t^8 + 1)*q^5 + O(q^7)
? [alpha, w] = params
% = [0, 8]
\end{verbatim}

This requires a few explanations: the {\tt mfslashexpansion} command asks
for $6$ terms of the Fourier expansion of $F|_k\ga$ with
$\ga=\psmm{1}{0}{2}{1}$, and the flag $1$ which follows asks
to give the result in algebraic form
if possible (set the flag to $0$ for complex floating point approximations). Thus, as for values of characters,
the variable ``{\tt t}'' which is printed is \emph{canonically}
$e^{2\pi i/16}$, root of $\Phi_{16}(t)=t^8+1$. The {\tt params} components
are $\al$ and $w$, and mean that the ``{\tt q}'' in the expansion should
be understood as $e^{2\pi i\tau/w}$ (so here $q=e^{2\pi i\tau/8}$),
and the expansion should be multiplied by $q^{\al}$ (here by $1$).

Here is a more complicated example which illustrates this:
\begin{verbatim}
? mf = mfinit([36,3,-4],0); F = mfbasis(mf)[1];
? Ser(mfslashexpansion(mf, F, g, 4, 1, &params), q)
% = Mod(-1/54*t^4 - 1/54*t, t^6 + t^3 + 1)
  + Mod(-1/3*t^5 - 1/3*t^2, t^6 + t^3 + 1)*q^2
  + Mod(-2/9*t^4, t^6 + t^3 + 1)*q^3 + O(q^5)
? [alpha, w] = params
% = [1/18, 9]
\end{verbatim}
Thus $w=9$ so $q=e^{2\pi i\tau/9}$, and $\al=1/18$ so the expansion must
be multiplied by $e^{2\pi i\tau/18}$. Note that the constant coefficient
$-t^4/54 - t/54 \in \Q(t)/(t^6 + t^3 + 1)$ in the expansion is
nonzero, so necessarily $\al>0$ otherwise $F$ would not be a cusp form.

Of course we can have ``raw'' expansions with approximate complex
coefficients:

\begin{verbatim}
? Ser(mfslashexpansion(mf, F, g, 4, 0), q)
% = (0.00321570699... - 0.01823718061... I)
  + (0.25534814770... - 0.21426253656... I)q^2 +...
\end{verbatim}
(Here we did not ask for {\tt params} since we already know it from the
previous computations.) Recall that all this is possible thanks to the
expression of $F$ as a linear combination of products of two Eisenstein
series. More generally for $\ga$ in $M_2^+(\Q)$, the result would be expressed
as $u^{k/2}f(\tau + v)$ and $f(q) = q^\alpha \sum_{m\ge0} a_m q^{m/w}$ as above
for some rational numbers $u$ and $v$ chosen so as to minimize the field of
definition of the $a_m$ (see the manual).

Atkin--Lehner operators are an important special case:
\begin{verbatim}
? mf = mfinit([32,4,8],0); Z = mfatkininit(mf,32);
? [mfB, M, C] = Z;
? M
% = 
[  1/8 -7/4]

[-1/16 -1/8]

? C
% = 0.35355339059327376220042218105242451964
\end{verbatim}
(here $C=8^{-1/2}$). This is a difficult case for the Atkin--Lehner operators
since first the level is not squarefree, and second the character, here
$8$ which represents $\lgs{8}{n}$, is not defined modulo $N/Q=32/32$.

The result involves a normalizing constant $C$ given above, essentially a
Gauss sum, such that the expansion of $C\cdot F|_k W_Q$ has the same field of
coefficients as $F$. A similar question for the more
general $F|_k\ga$ is answered by Theorem \ref{thmfga}, but the given value
may not be optimal. In general, $C\cdot F|_k W_Q$ belongs to a different
space than $F$ (the Nebentypus becomes $\overline{\chi_Q} \chi_{N/Q}$ in
integral weight, and a similar formula holds in non-integral weight). The matrix $M$ expresses $C \cdot F_i |_k W_Q$ when
$(F_i)$ is a basis of {\tt mf} in terms of a basis of that other space
{\tt mfB}. The operator can also be applied to an individual form:
\begin{verbatim}
? F = mfbasis(mf)[1]; G = mfatkin(Z,F);
? mfser(G,10)
% = 1/4*q + 7/2*q^3 + 7*q^5 + 2*q^7 - 1/4*q^9 + ...
\end{verbatim}
This returns $G = C\cdot F|_k W_{32}$.

\subsection{Numerical Applications}
  
It is now easy to compute \emph{period polynomials}
$$P(F,X)=\int_0^{i\infty}(X-\tau)^{k-2}F(\tau)\,d\tau$$
for modular forms of integral weight $k\ge2$.
Indeed, in addition to the Fourier expansion at infinity, this only requires
the expansion of $F|_kW_N$ with $W_N=\psmm{0}{-1}{N}{0}$.
In particular we can compute in complete generality \emph{special values}
and \emph{periods}.

More generally, we can compute numerically \emph{general} period polynomials
and \emph{modular symbols}
$$\int_{s_1}^{s_2}(X-\tau)^{k-2}(F|_k\ga)(\tau)\,d\tau$$
for any cusps $s_1$ and $s_2$ and any $\ga\in M_2^+(\Q)$, or even for
any $s_1$ and $s_2$ in the completed upper-half plane $\H$.

Thanks to these symbols, we can also compute general \emph{Petersson products}
of modular forms of integral weight $k\ge2$ by using Haberland-type formulas
such as Theorem \ref{thmpet}:

\begin{verbatim}
? mf = mfinit([11,2],0); F = mfbasis(mf)[1];
? FS = mfsymbol(mf,F); mfsymboleval(FS,[0,oo])
% = 0.040400186918863279214419198537327720301*I
? mfsymboleval(FS, [1/2,2*I])
% = -0.16160130270129178714056342927215612575*I
? mfpetersson(FS) \\ <F,F>
% = 0.0039083456561245989852473854813821138618
\end{verbatim}

\begin{verbatim}
? mf = mfinit([23,2],1);
? BS = [mfsymbol(mf,f) | f <- mfbasis(mf)];
? [mfpetersson(f,g) | f<-BS; g<-BS]
% = [0.0095931508727672866790131897867245345540,
     -0.0066920429957575620051313153184106231192,
     -0.0066920429957575620051313153184106231192,
      0.016285193868524848684144505105135157673]
\end{verbatim}
We can also evaluate a form near the real axis: for forms over the full
modular group $\G$ or small index, we can use modular transformations
\emph{in the group} to significantly increase imaginary parts.
In general, this is not possible, but now it is easy since we can always
use the whole of $\G$: in fact, we can always reduce to
$\Im(\tau)\ge 1/(2N)$, which is almost optimal:

\begin{verbatim}
? \p57
? mf = mfinit([12,4],1); F = mfbasis(mf)[1];
? ev(m) = mfeval(mf, F, 1/Pi+I/10^m); 
? ev(6)
% = -89811.0493... -58409.9409...*I
? ev(7)
% = 4.821... E-52 + 6.788... E-52*I
? ev(8)
% = -1.79763... E-69 + 2.6450... E-69*I
? ev(9)
% = 357873461.23... - 264528426.36...*I
? ev(10)
% = 0.3966...E18 - 1.6429...E18*I
\end{verbatim}

  Note that $|ev(m)|$ seems to tend to infinity with $m$, but with a pronounced
  ``dip'' around $m=7$ and $m=8$. This is not specific to the modular form
  $F$ but probably to the diophantine approximation properties of $1/\pi$,
  and in particular to its very close convergent $113/355$.
 

\end{document}